\documentclass[10pt]{article} 
\usepackage{amsmath,amssymb}
\usepackage{booktabs} 
\usepackage{array} 
\usepackage{paralist} 
\usepackage{verbatim} 
\usepackage{subfig} 

\usepackage{fancyhdr} 
\usepackage{sectsty}
\allsectionsfont{\sffamily\mdseries\upshape} 

\usepackage[nottoc,notlof,notlot]{tocbibind} 
\usepackage[titles,subfigure]{tocloft} 

\title{On some cryptological schemes for r-person secret ballot and r-person authentication}
\author{BongJu Kim}
\date{} 

\begin{document}
\maketitle

\section{A cryptological schemes using the partition identities}
\subsection{A membership verification using partition identities}

Suppose that $r$ members want to keep their fortune in a bank.  Approaching to their fortune requires all member's permission and it will be happened many times. In the membership verification processes, they don't believe the bank employee's decision whether  each claimant is really one of the members and don't believe each other who claims his/her membership and worry about that in the communication processes with the bank, informations they exchange may leak.

We introduce a scheme for the membership verification using a kind of partition identities. In $\cite{2}$, Kim introduced some kind of partition number identities. A feature of this identities is that although we change the base set(or the set of parts)\footnote{It is a subset of the set of natural numbers. When we consider partitions on a base set(set of parts), we only use elements of the base set.}, this identities always true. We introduce one of the identities briefly in Section 2 and suggest cryptological schemes, can be used to r-person authentication and r-person secret ballot, using the identities in Section 3. The schemes based on the difficulty of computing $u$, $v$ from the values of $p^{A}_{\alpha}(u)p^{A}_{\alpha}(v)$ for some base sets $A$. For brief explanation, see the second chapter.

Here is a scheme for the membership verification.

Setting. For $r$ members, They put a valuable in a safe and share the key. Then they keep it in a safe-deposit box at a bank. Assume that the safe is really safe.

Step 1. A bank employee chooses natural numbers $n_1$, $n_2$ and $\alpha$ such that $gcd(n_1, \alpha)=gcd(n_2, \alpha)=1$.\footnote{If $\alpha | n$, the solution matrix contains a divisor of $n$.} Then he/she computes the solution matrix $(a^{n_\mu,\alpha}_{ij})$ for $n_1$, $n_2$ and $\alpha$ so that he/she gets the partition identity
\[p^{A}(n_\mu)=\sum_{i\geq 1}\prod_{j\geq 1}p^{A}_{\alpha}(a^{n_\mu,\alpha}_{ij}).\]Since the solution matrix $(a^{n,\alpha}_{ij})$ contains $n$, we do not want to reveal it.
The bank employee computes the product
\[(p^{A}(n_1)-p^{A}_{\alpha}(n_1))(p^{A}(n_2)-p^{A}_{\alpha}(n_2))\]
and expend it. The expansion is composed of the sum of 
\[\prod_{j\geq 1}p^{A}_{\alpha}(a^{n_1,\alpha}_{i_tj})\prod_{j\geq 1}p^{A}_{\alpha}(a^{n_2,\alpha}_{i_sj}).\]
We rewrite the expansion by
\[(p^{A}(n_1)-p^{A}_{\alpha}(n_1))(p^{A}(n_2)-p^{A}_{\alpha}(n_2))=\sum_{i\geq 1}\prod_{j\geq 1}p^{A}_{\alpha}(b_{ij}).\]
The bank employee keeps secret $n_1$, $n_2$ and the identities. He/She let the $r$-members know $\alpha$.

Step 2. The bank employee chooses some parts of terms of 
\[\sum_{i\geq 1}\prod_{j\geq 1}p^{A}_{\alpha}(b_{ij})\]
 so that any $b_{ij}$ in a chosen term is far enough from $1$ and $n_1, n_2$. The bank employee properly distributes among the $r$ members in the form of a product of about half of original terms\footnote{Revealing all term of a product can be helpful to compute $n_1$ or $n_2$. Revealing a single partition number $p^{A}_{\alpha}(b_{ij})$ also should be avoided.} such as $\prod_{k=1}^{K_1}p^{A}_{\alpha}(b_{i_1j_k})$ for the first member, $\cdots$, $\prod_{k=1}^{K_r}p^{A}_{\alpha}(b_{i_1j_k})$ for the $r$th member\footnote{More complicated form of distribution such as the bank employee gives an ordered pair $(\prod_{k=1}^{K_1}p^{A}_{\alpha}(b_{i_1j_k}), \prod_{k=1}^{K_2}p^{A}_{\alpha}(b_{i_2j_k}), \cdots, \prod_{k=1}^{K_t}p^{A}_{\alpha}(b_{i_tj_k}))$ for a member also possible as long as any single term of distributed ordered pair is a product of at least two terms of identity and preserving the original from i.e. a member receives a product $p^{A}_{\alpha}(u)p^{A}_{\alpha}(v)$ if and only if there exists $i$ such that $u$, $v\in \{b_{ij}\;|\; j\;\}$.}. Note that it is impossible that restoring $n_1$, $n_2$ by only a part of the identity. This process should be performed in safe situation without any possibility of eavesdropping and the members should not know each other's terms.

Step 3. Every three years, The bank employee and the claimants chooses randomly a large base set $A^*$. The bank employee computes \[(p^{A^*}(n_1)-p^{A^*}_{\alpha}(n_1))(p^{A^*}(n_2)-p^{A^*}_{\alpha}(n_2))\] and the terms which hadn't been distributed in the base set $A^*$. Each claimant computes the combination of partition numbers, which they received, in the base set $A^*$. Each claimant shares his/her result. 

Step 4. The bank employee informs of the values of the undistributed terms and the method to combine each claimant's results with the values. The claimants combines the partition numbers from the persons according to the method. The bank employee informs of $(p^{A^*}(n_1)-p^{A^*}_{\alpha}(n_1))(p^{A^*}(n_2)-p^{A^*}_{\alpha}(n_2))$ and the claimants informs of their final result. They check them. 

Step 5. Repeat Step 3 and Step 4 for several times. 

Step 6. If the all tests for several base sets passed, The bank employee opens the safe-deposit box. If they fail to the authentication process, The bank employee reveals the values of distributed terms. The claimants find who's process was fail and expel him/her/them.

In Step 3 and Step 4, each claimant knows only $(p^{A^*}(n_1)-p^{A^*}_{\alpha}(n_1))(p^{A^*}(n_2)-p^{A^*}_{\alpha}(n_2))$ and the other's results for several base sets so he/she can not compute the other's $b_{ij}$ nor $n_1$, $n_2$. Therefore each claimant can not perform himself/herself nor make a person can perform as if he/she is a member except a renouncement of his/her membership. An eavesdropper of this process can not perform as if he/she is one of the member since the eavesdropper can not know more information than a claimant. Therefore, in the situation that the valuable, stored in a system which has same function with The bank employee in the scheme, is an digital information and encrypted the member's key, all communications made in on-line, this scheme is safe from an eavesdropper or fake claimant.

\subsection{A scheme for a secret ballot}
Let us think about a secret ballot in a community. The community have $r$ decision-makers who elected by all members of the community. There is also an election commission for a secret ballot which does not contains a decision-maker\footnote{In a representative democracy, for example, all citizens elect members of the National Assembly. In the National Assembly, all members of the National Assembly elect the Chairman of the National Assembly(election commission for decision making) and a standing committee(decision-makers). A bill should be passed a standing committee before the Assembly plenary session. The Chairman of the National Assembly can not participate in a standing committee's decision making.}. The election commission gives a suffrage to each decision-maker. In a secret ballot, the members of the community want to know the number of ayes and the number of nays, but it should not be known that who is the nay or the aye except himself/herself.

In this situation, the scheme in Section 3.1 would be useful after some modification. Here is a scheme for a secret ballot.

Setting. For decision-makers of $r$ members, the decision-makers meet in person and does not reveal any information unless specifically noted. The election commission does not contains a decision-maker and also does not reveal any information unless specifically noted. The communication between the election commission and the decision-makers is opened to the public. The number $r$ of the decision-makers is also known to the public.

Step 1. The election commission chooses natural numbers $n_1$, $n_2$ and $\alpha$ such that $gcd(n_1, \alpha)=gcd(n_2, \alpha)=1$.\footnote{If $\alpha | n$, the identity contains a divisor of $n$.} Then he/she computes the solution matrix $(a^{n_\mu,\alpha}_{ij})$ for $n_1$, $n_2$ and $\alpha$ so that he/she gets the partition identity
\[p^{A}(n_\mu)=\sum_{i\geq 1}\prod_{j\geq 1}p^{A}_{\alpha}(a^{n_\mu,\alpha}_{ij}).\]Since the solution matrix $(a^{n,\alpha}_{ij})$ contains $n$, we do not want to reveal it.
The election commission computes the product
\[(p^{A}(n_1)-p^{A}_{\alpha}(n_1))(p^{A}(n_2)-p^{A}_{\alpha}(n_2))\]
and expend it. The expansion is composed of the sum of 
\[\prod_{j\geq 1}p^{A}_{\alpha}(a^{n_1,\alpha}_{i_tj})\prod_{j\geq 1}p^{A}_{\alpha}(a^{n_2,\alpha}_{i_sj}).\]
We rewrite the expansion by
\[(p^{A}(n_1)-p^{A}_{\alpha}(n_1))(p^{A}(n_2)-p^{A}_{\alpha}(n_2))=\sum_{i\geq 1}\prod_{j\geq 1}p^{A}_{\alpha}(b_{ij}).\]
The election commission keeps secret $n_1$, $n_2$ and the identities. He/She let the public know $\alpha$.

Step 2. The election commission chooses two sets \[ E:=\{i_1, i_2, \cdots, i_r, i_{r+1}, i_{r+2}, \cdots, i_{2r}\} \] and \[B':=\{b'_{ij}\}\] from $B:=\{b_{ij}\}$ where any $b'_{kj}$ in the products is far enough from $1$ and $n_1, n_2$. We define for $k\in E$,
\[v'_k(A):=\prod_{j\geq 1}p^{A}_{\alpha}(b'_{kj}),\;\; b'_{kj}\in B',\] \[ v_k(A):=\prod_{j\geq 1}p^{A}_{\alpha}(b_{kj}),\;\; b_{kj}\in B'^c\] and for $k\in E^c$,\[u'_k(A):=\prod_{j\geq 1}p^{A}_{\alpha}(b'_{kj}),\;\; b'_{kj}\in B',\] \[ u_k(A):=\prod_{j\geq 1}p^{A}_{\alpha}(b_{kj}),\;\; b_{kj}\in B'^c.\]
 Note that \[(p^{A}(n_1)-p^{A}_{\alpha}(n_1))(p^{A}(n_2)-p^{A}_{\alpha}(n_2))\] can be expressed by\[\sum_{k\in E}v'_k(A)v_k(A)+\sum_{k\in E^c}u'_k(A)u_k(A).\] The election commission distributes among the $r$ decision-makers such that the $s$th decision-maker has two term\footnote{The number of terms a decision-maker has can be grater by choosing the set $E$ properly.}\[v'_{i_s}, v'_{i_{r+s}}\] and opens \[\{u'_k(A),\;|\;k\in E^c\}\] to the public.\footnote{This information can be helpful in a inspection without anonymity and this may make a community member who is not a member of the election commission nor the decision-makers more feels that he/she has something to do with the vote.} Note that it is impossible that restoring $n_1$, $n_2$ by $u'_k, v'_k$ only a part of the identity. This process should be performed in safe situation without any possibility of eavesdropping and the decision-makers should not know each other's terms.

Step 3. In a secret ballot, the election commission and the decision-makers chooses randomly a large base set $A^*$ so that the decision-maker's own computation is not $0$. If a decision-maker's computation is $0$, the decision-maker can not express a disagreement with anonymity. They inform $A^*$ of the public. The election commission computes \[(p^{A^*}(n_1)-p^{A^*}_{\alpha}(n_1))(p^{A^*}(n_2)-p^{A^*}_{\alpha}(n_2))\] and two ordered pairs of natural numbers\[(v_{i_s}(A^*), v_{i_{r+s}}(A^*))_{s=1}^r, (u_{k}(A^*))_{k\in E^c}\]and informs the public of the ordered pairs. Each decision-maker computes \[V_s(A^*):=v'_{i_s}(A^*)v_{i_s}(A^*)+v'_{i_{r+s}}(A^*)v_{i_{r+s}}(A^*)\]with the values, which are announced by the election commission, of the undistributed terms.  If a decision-maker wants to express his/her agreement, he/she adds $1$ to his/her original computation $V_s$. Each aye informs the decision-makers of his/her \[W_s(A^*):=V_s(A^*)+1\] and each nays informs of the decision-makers of his/her \[W_s(A^*):=V_s(A^*).\] They share their results among only the decision-makers.

Step 4. The decision-makers compute \[\sum_{s=1}^rW_s(A^*)+\sum_{k\in E^c}u'_k(A)u_k(A).\] The election commission reveals $(p^{A^*}(n_1)-p^{A^*}_{\alpha}(n_1))(p^{A^*}(n_2)-p^{A^*}_{\alpha}(n_2))$ to the public and, at the same time, the decision-makers reveal their final result  to the public. They compute the number of the ayes \[y=\sum_{s=1}^rW_s(A^*)+\sum_{k\in E^c}u'_k(A)u_k(A)-(p^{A^*}(n_1)-p^{A^*}_{\alpha}(n_1))(p^{A^*}(n_2)-p^{A^*}_{\alpha}(n_2))\] and the number of nays\[r-y.\]

Step 5. If there are more than two items in the vote, Repeat Step 3 and Step 4 for each item.

Inspection with anonymity. A decision-maker can have a fraudulent vote by adding a number other than $1$ or $0$ such as  a positive integer grater than $1$ or a negative integer in Step 3. The community wants to inspect whether a decision-maker had had a fraudulent vote and also wants to keep it's anonymity. The number of the case of fair vote is $2^r$ since the number of the decision-makers is $r$ and a decision-maker chooses $1$ or $0$. They determine a hash function\[h:\mathbb{Z}^r \longrightarrow \mathbb{Z}.\]A fair vote is 
\[ (V_s(A^*))_{s=1}^r+w,\]where $w\in\{0, 1\}^r$. The election commission makes a table of hashes of all possible $2^r$ cases of fair votes\footnote{It is reasonable if $r$ is not too large. In Korea now(2016), the number of members of a standing committee(decision-makers) does not exceed 31 and the average is about 22 except Special Committee on Budget and Accounts(50 members).} and informs the public of the table\[\{h((V_s(A^*))_{s=1}^r+w) \;|\; w\in\{0, 1\}^r\}\subset\mathbb{Z}.\]  
Then the decision-makers compute the hash of the vote \[h((W_s(A^*))_{s=1}^r)\] and inform the public of this value. The community checks if one of the values of the table is equals to this value. If there is no match, then a decision-maker had added a number other than $0$, $1$, i.e. it was a fraudulent vote.

\subsection{A scheme for the unanimity rule which can hide the number of voters}
Let us think about the community, mentioned at the previous section, wants a decision making of the decision-makers according to the unanimity rule. The election commission gives a suffrage to each decision-maker. Every time the community needs a decision making according to the unanimity rule, they hold a secret ballot. In a secret ballot, the members of the community want to know whether there is an objector, but it should not be known that who is the objector except the objector himself/herself.

Here is a scheme for the community. This scheme do not reveal the number of the objectors and the number of the decision-makers. In addition, it is possible that the objectors reveal and prove their disagreement without a help of the election commission or the other decision-makers if they want to do that. 

Setting. For decision-makers of $r$ members, the decision-makers meet in person and does not reveal any information unless specifically noted. the election commission does not contains a decision-maker and also does not reveal any information unless specifically noted. The communication between the election commission and the decision-makers is opened to the public. They do not reveal the number $r$ of the decision-makers. 

Step 1. Same with Step 1 in Section 1.2.

Step 2. Same with Step 2 in Section 1.2.

Step 3. In a secret ballot, the election commission and the decision-makers chooses randomly a large base set $A^*$. The election commission computes \[(p^{A^*}(n_1)-p^{A^*}_{\alpha}(n_1))(p^{A^*}(n_2)-p^{A^*}_{\alpha}(n_2))\]and the undistributed terms in the base set $A^*$. Each decision-maker computes the combination of partition numbers, which they received, in the base set $A^*$. If a decision-maker want to express his/her disagreement, he/she changes randomly his/her result so that it is different from before. Each decision-maker shares his/her result among the decision-makers.

Step 4. The election commission informs the public of the values of unrevealed terms. The decision-makers combine the shared values with the values of unrevealed terms. The election commission reveals $(p^{A^*}(n_1)-p^{A^*}_{\alpha}(n_1))(p^{A^*}(n_2)-p^{A^*}_{\alpha}(n_2))$ to the public and, at the same time, the decision-makers reveal their final result to the public. They check if the two values are equal.

Step 5. Repeat Step 3 and Step 4 for several times. 

Step 6. If the all tests for several base sets passed, the community judges that the decision-makers made an unanimous agreement in this vote. If they failed in a process, the community judges that a decision-maker had disagreed in this vote.\footnote{If two or more decision-maker changes randomly his/her result, it is possible that the final value still remains same with $(p^{A^*}(n_1)-p^{A^*}_{\alpha}(n_1))(p^{A^*}(n_2)-p^{A^*}_{\alpha}(n_2))$.}

If one wants to know who was disagreed, he/she should know both each decision-maker's result and the terms each decision-maker received in Step 2. It is possible only when one of the election commission or the decision-makers leaked his/her secret so that someone knows the informations from both of them.

If decision-makers who had disagreed in a vote want to reveal their disagreements to the public, they gather all objectors in the vote and reveal each objector's original value and the other decision-maker's values they computed in Step 3. Then the public can compute the final value and check if it equals to $(p^{A^*}(n_1)-p^{A^*}_{\alpha}(n_1))(p^{A^*}(n_2)-p^{A^*}_{\alpha}(n_2))$. It proves their disagreements. It may be useful if the objectors encounter some disadvantages due to their disagreements. After the disagreement of them is proved, a judicial authority of the community investigates the election commission and the other decision-makers.
 
\section{Partition number Identities}
\subsection{A partition number} 
Let $A$ be a subset of the set of natural numbers. We define a partition number \[p^A_\alpha(n)\] by the number of partitions of $n$ into elements of $A$ such that the number of equal parts is less than or equals to $\alpha$. If there is no restriction of the number of equal parts, We denote by \[p^A(n).\]

\subsection{A solution set} Let $n$ and $\alpha$ be natural numbers. We consider all representations of $n$ by combinations of powers of the $\alpha$, i.e. all solutions of the equation \[n=\sum_{i\geq1}N_i\alpha^i.\] For example, $n$ can be represented as following
 \[n=a^{n}_{11}+2a^{n}_{12}+4a^{n}_{13}+\cdots\\
    =a^{n}_{21}+2a^{n}_{22}+4a^{n}_{23}+\cdots\\
    =\cdots. \]
Now, we gather all solutions and make the solution matrix \[(a^{n,\alpha}_{ij}).\]
\subsection{Partition number identities}
Let $A$ be a subset of the set of natural numbers. Let $(a_{ij}^{n,\alpha})$ be the solution matrix of $n=\sum_{i\geq 0}(\alpha +1)^{i}N_{i}$.  Let $n$ and $\alpha$ be natural numbers. Then
\[
p^{A}(n)=\sum_{i\geq 1}\prod_{j\geq 1}p^{A}_{\alpha}(a^{n,\alpha}_{ij})
\]
for all $A$.
\subsection{Examples}
We show that an identity for $n=10$, $\alpha=1$ holds for three different base sets
\[\begin{array}{l}
 p^{A}(10) \\
 =p^{A}_{1}(1)p^{A}_{1}(1)+p^{A}_{1}(1)p^{A}_{1}(2)+p^{A}_{1}(1)p^{A}_{1}(3) \\
  +p^{A}_{1}(5)+p^{A}_{1}(2)p^{A}_{1}(2)p^{A}_{1}(1)\\ 
  +p^{A}_{1}(2)p^{A}_{1}(4)+p^{A}_{1}(2)p^{A}_{1}(1) \\ 
 +p^{A}_{1}(2)p^{A}_{1}(2)+p^{A}_{1}(4)p^{A}_{1}(1)p^{A}_{1}(1) \\ 
 +p^{A}_{1}(4)p^{A}_{1}(3)+p^{A}_{1}(6)p^{A}_{1}(2)+p^{A}_{1}(6)p^{A}_{1}(1) \\ 
 +p^{A}_{1}(8)p^{A}_{1}(1)+p^{A}_{1}(10).
 \end{array} 
 \]
 Now, we check the identity for three sets of parts.\\1) $A=\{p\;|$ $p$ is a prime$\}.$\\Let us calculate partition numbers.\\
 $
 \begin{array}{ll}
 10= & 5+5 \\ 
  & 3+2+5 \\ 
  & 3+3+2+2 \\ 
  & 7+3 \\ 
  & 2+2+2+2+2
 \end{array} 
 $, $
 \begin{array}{ll}
 p^{A}_{1}(1)=0 & p^{A}_{1}(5)=2 \\ 
 p^{A}_{1}(2)=1 & p^{A}_{1}(6)=0 \\ 
 p^{A}_{1}(3)=1 & p^{A}_{1}(8)=1 \\ 
 p^{A}_{1}(4)=0 & p^{A}_{1}(10)=2.
 \end{array} 
 $\\

So, $p^{A}(10)=5$. Since $p^{A}_{1}(1)=p^{A}_{1}(4)=p^{A}_{1}(6)=0$, 
 \[
 p^{A}(10)=p^{A}_{1}(5)+p^{A}_{1}(2)p^{A}_{1}(2)+p^{A}_{1}(10)=2+1+2=5.
 \]
\\
2) $A=\{n^2\} $.\\Let us calculate the partition numbers.\\
$
 \begin{array}{ll}
 10= & 1+9 \\ 
  & 1+1+1+1+1+1+1+1+1+1 \\ 
  & 4+1+1+1+1+1+1 \\ 
  & 4+4+1+1
 \end{array} 
$, $
 \begin{array}{ll}
 p^{A}_{1}(1)=1 & p^{A}_{1}(5)=1 \\ 
 p^{A}_{1}(2)=0 & p^{A}_{1}(6)=0 \\ 
 p^{A}_{1}(3)=0 & p^{A}_{1}(8)=0 \\ 
 p^{A}_{1}(4)=1 & p^{A}_{1}(10)=1.
 \end{array} $

So, $p^{A}(10)=4$. Since $p^{A}_{1}(2)=p^{A}_{1}(3)=p^{A}_{1}(6)=p^{A}_{1}(8)=0$,
\[\begin{array}{ll}
p^{A}(10)&=p^{A}_{1}(1)p^{A}_{1}(1)+p^{A}_{1}(5)+p^{A}_{1}(4)p^{A}_{1}(1)p^{A}_{1}(1)+p^{A}_{1}(10)\\&=1+1+1+1=4.
\end{array}\]
\\
3) $A=\{n\;|$ $n$ is an odd number$\}$.\\Let us calculate the partition numbers.\\
$
 \begin{array}{ll}
 10= & 1+1+1+1+1+1+1+1+1+1 \\ 
  & 3+1+1+1+1+1+1+1 \\ 
  & 5+1+1+1+1+1 \\
  & 7+1+1+1 \\
  & 9+1 \\
  & 3+3+1+1+1+1 \\
  & 3+3+3+1 \\
  & 5+5 \\
  & 7+3 \\
  & 3+5+1+1
 \end{array} 
$, $
 \begin{array}{ll}
 p^{A}_{1}(1)=1 & p^{A}_{1}(5)=1 \\ 
 p^{A}_{1}(2)=0 & p^{A}_{1}(6)=1 \\ 
 p^{A}_{1}(3)=1 & p^{A}_{1}(8)=2 \\ 
 p^{A}_{1}(4)=1 & p^{A}_{1}(10)=2.
 \end{array}$\\So, $p^{A}(10)=10$. Since $p^{A}_{1}(2)=0$,\\
 
$\begin{array}{ll}
  p^{A}(10) & =p^{A}_{1}(1)p^{A}_{1}(1)+p^{A}_{1}(1)p^{A}_{1}(3)+p^{A}_{1}(5) \\
  & \,\,+p^{A}_{1}(4)p^{A}_{1}(1)p^{A}_{1}(1)+p^{A}_{1}(4)p^{A}_{1}(3) \\
  & \,\,+p^{A}_{1}(6)p^{A}_{1}(1)+p^{A}_{1}(8)p^{A}_{1}(1)+p^{A}_{1}(10) \\
  & \,\,=1+1+1+1+1+1+2+2 \\
  & \,\,=10. 
 \end{array}$

\section{Closing remarks}
A r-person authentication can be performed by using an existing authentication scheme which based on difficult problem such as integer factoring or discrete logarithm or based on elliptic curves. But these problems studied for many years and some effective methods for the problems are known. For example, Pollard-Strassen prime factorization$\cite{1}$, method Especially, it is known that there are polynomial-time algorithms for a quantum computer solving integer factoring and discrete logarithm$\cite{3}$. However, it seems difficult that one compute $u$, $v$ from the values of $p^{A^*}_{\alpha}(u)p^{A^*}_{\alpha}(v)$ for some base sets $A^*$.

To use the schemes in practice, the method to choose proper $n_1$, $n_2$, $\alpha$ and $A^*$ should be invented in consideration of computing ability at that moment.


\begin{thebibliography}{widestlabel}

\bibitem{1} Hardy, K.; Muskat, J. B.; and Williams, K. S. {\it A Deterministic Algorithm for Solving $n=fu^2+gv^2$ in Coprime Integers u and v}. Math. Comput. 55, 327-343, 1990

\bibitem{2} Kim, BongJu {\it On some elementary partition number identities}.  arXiv:1803.08095, 2012

\bibitem{3} Shor, Peter W. {\it Polynomial-time algorithms for prime factorization and discrete logarithms on a quantum computer}. SIAM J. Comput. 26 (1997), no. 5, 1484–1509.

\end{thebibliography}
\end{document}